\documentclass[12pt,reqno]{amsart}
\usepackage{amsmath}
\usepackage{amsxtra}
\usepackage{amscd}
\usepackage{amsthm}
\usepackage{amsfonts}
\usepackage{amssymb}
\usepackage{eucal}

\theoremstyle{definition}

\theoremstyle{remark}

\theoremstyle{remark}

\numberwithin{equation}{section}

\let\Maketitle\maketitle
\def\maketitle{\hrule height0pt
\Maketitle\thispagestyle{empty}\let\maketitle\empty}

\let\geq\geqslant

\def\>{\relax\ifmmode\mskip.666667\thinmuskip\relax\else\kern.111111em\fi}
\def\<{\relax\ifmmode\mskip-.333333\thinmuskip\relax\else\kern-.0555556em\fi}

\def\vv#1>{\vadjust{\vskip#1\baselineskip}}
\def\vvn#1>{\vadjust{\nobreak\vskip#1\baselineskip\nobreak}}
\def\vvgood{\vadjust{\penalty-500}}

\def\XXX/{{\sl XXX\/}}

\let\nc\newcommand

\nc{\wt}{\widetilde}
\nc{\on}{\operatorname}
\nc{\ch}{\mbox{ch}}
\nc{\Z}{\mathbb Z}
\nc{\C}{\mathbb C}
\nc{\R}{\mathbb R}
\nc{\pone}{\mathbb{CP}^1}
\nc{\p}{\partial}
\nc{\F}{\mathcal F}
\nc{\arr}{\rightarrow}
\nc{\larr}{\longrightarrow}
\nc{\al}{\alpha}
\nc{\ri}{\rangle}
\nc{\lef}{\langle}
\nc{\W}{\mathcal W}
\nc{\la}{\lambda}
\nc{\ep}{\epsilon}
\nc{\eps}{\varepsilon}
\nc{\Om}{\Omega}
\nc{\om}{\omega}
\nc{\La}{\Lambda}
\nc{\Ga}{\Gamma}
\nc{\ga}{\gamma}

\nc{\D}{\mathcal D}
\nc{\A}{{\mathbb A}}
\nc{\g}{\mathfrak g}
\nc{\h}{\mathfrak h}
\nc{\mg}{\mathfrak m}
\nc{\n}{\mathfrak n}
\nc{\z}{\mathfrak Z}
\nc{\mc}{\mathcal}
\nc{\M}{\mathcal M}
\nc{\N}{\widehat\n}
\nc{\G}{\widehat\g}
\nc{\De}{\Delta_+}
\nc{\gt}{\widetilde{\g}}
\nc{\one}{\mathbf 1}
\nc{\qp}{q^{\frac{k}{2}}}
\nc{\qm}{q^{-\frac{k}{2}}}
\nc{\qn}{\frac{[m]_q^2}{[2m]_q}}
\nc{\Res}{\on{Res}}

\nc{\Wr}{{\rm Wr}}

\nc{\bean}{\begin{eqnarray}}
\nc{\eean}{\end{eqnarray}}
\nc{\be}{\begin{equation*}}
\nc{\ee}{\end{equation*}}
\nc{\bea}{\begin{eqnarray*}}
\nc{\eea}{\end{eqnarray*}}
\nc{\beq}{\begin{equation}}
\nc{\eeq}{\end{equation}}
\nc{\bs}{\boldsymbol}
\nc{\Ref}[1]{{$($\ref{#1}$)$}}

\nc{\glN}{\mathfrak{gl}_N}
\nc{\glt}{\mathfrak{gl}_2}

\nc{\End}{{\rm End}}
\nc{\T}{{\bs p}}
\nc{\nash}{M_{\bs \La}[\La^{(\infty)}]}
\nc{\Sing}{{\rm Sing\,}}
\nc{\snash}{\Sing M_{\bs \La}[\La^{(\infty)}]}
\nc{\Snash}{\Sing L_{\bs \La}[\La^{(\infty)}]}
\nc{\sh}{{\sigma}}

\nc{\AD}{A_{\T,D}}
\nc{\AT}{A_{\T,T}}
\nc{\Dh}{\mathcal{D}_{\bs h}}
\nc{\Dhp}{\mathcal{D}_{\bs h(\bs p)}}

\nc{\qa}{p(x,\bs a)}
\nc{\qat}{\tilde {p}(x,\tilde {\bs a})}

\nc{\qap}{p(x,\bs a(\bs p))}

\nc{\Ht}{{\widetilde{H} }}

\begin{document}

\title[ Bethe algebra and algebra of functions]
{Bethe algebra and algebra of functions on the space\\[3pt]
of differential operators of order two\\[3pt]
with polynomial solutions }

\author[E.\,Mukhin, V.\,Tarasov, and A.\>Varchenko]
{E.\,Mukhin ${}^{*,1}$, \,V.\,Tarasov ${}^{*,\star,2}$,
\,and \,A.\>Varchenko {${}^{**,3}$} }

\thanks{${}^1$\ Supported in part by NSF grant DMS-0601005}
\thanks{${}^2$\ Supported in part by RFFI grant 05-01-00922}
\thanks{${}^3$\ Supported in part by NSF grant DMS-0555327}

\maketitle

\centerline{\it ${}^*$Department of Mathematical Sciences,
Indiana University -- Purdue University,}
\centerline{\it Indianapolis, 402 North Blackford St, Indianapolis,
IN 46202-3216, USA}
\smallskip
\centerline{\it $^\star$St.\,Petersburg Branch of Steklov Mathematical
Institute}
\centerline{\it Fontanka 27, St.\,Petersburg, 191023, Russia}
\smallskip
\centerline{\it ${}^{**}$Department of Mathematics, University of
North Carolina at Chapel Hill,} \centerline{\it Chapel Hill, NC
27599-3250, USA} \medskip

\medskip

\begin{abstract}
We show that the following two algebras are isomorphic.
The first is the algebra $A_P$ of functions on the scheme of monic linear
second-order differential operators on $\C$ with prescribed regular
singular points at $z_1,\dots, z_n, \infty$, prescribed exponents
$\La^{(1)}, \dots, \La^{(n)}, \La^{(\infty)}$ at the singular points,
and having the kernel consisting of polynomials only.
The second is
the Bethe algebra of commuting linear operators, acting on the
vector space $\Sing\, L_{\La^{(1)}} \otimes \dots \otimes
L_{\La^{(n)}}[\La^{(\infty)}]$ of singular vectors of weight
$\La^{(\infty)}$ in the tensor product of finite dimensional
polynomial $\glt$-modules with highest weights $\La^{(1)},
\dots, \La^{(n)}$.
\end{abstract}

\section{Introduction}
\subsection{}
There is a classical connection between Schubert calculus and representation
theory of the Lie algebra $\glN$. Let $V$ be a vector space. Then Schubert
cycles in the Grassmannian of $N$-dimensional subspaces of $V$ are labeled
by highest weights of polynomial irreducible $\glN$-modules and if the
intersection of several cycles is finite, then the intersection number is equal
to the multiplicity of the unique one-dimensional representation in the tensor
product of the corresponding polynomial finite-dimensional $\glN$-modules.
It is a challenge to understand in a deeper way this numerological relation,
see~\cite{F}, \cite{B}.

In this paper we prove a result which may help to comprehend better
the interrelation of Schubert calculus and representation
theory. Namely, for $N=2$ under certain conditions, we identify the
algebra of functions on the intersection of Schubert cycles with the
Bethe algebra of linear operators acting on the multiplicity space
of the one-dimensional subrepresentation.

\subsection{}
Let $\La^{(1)}, \dots, \La^{(n)}, \La^{(\infty)}$ be dominant integral
$\glN$-weights. Consider
the tensor product $L_{\bs \La}=L_{\La^{(1)}}\otimes\dots\otimes L_{\La^{(n)}}$
of $n$ polynomial irreducible finite-dimensional $\glN$-modules with highest
weights $\La^{(1)}, \dots, \La^{(n)}$, respectively.
Let $\Sing\,L_{\bs \La}[\La^{(\infty)}] \subset L_{\bs \La}$
be the subspace of singular vectors of weight $\La^{(\infty)}$.
Fix $n$ distinct complex numbers $z_1,\dots,z_n$.
Then the theory of the integrable Gaudin model provides us with a collection of
commuting linear operators on that space, the operators being called the higher
Gaudin Hamiltonians or the higher transfer matrices. The unital algebra 
$A_L$ of
endomorphisms of $\Sing\,L_{\bs \La}[\La^{(\infty)}]$, generated by
the higher Gaudin Hamiltonians, is called the Bethe algebra.

\medskip

Thus, given a set of $n+1$ highest weights $\La^{(1)},
\dots, \La^{(n)},\La^{(\infty)}$ and a collection of complex numbers
$z_1, \dots, z_n$ we construct the vector space $\Sing\,L_{\bs
\La}[\La^{(\infty)}]$ and the commutative Bethe
algebra of linear operators acting on that space.

\medskip
There is another construction which starts with the same initial data.
Having a set of highest weights $\La^{(1)}, \dots,\La^{(n)},\La^{(\infty)}$
as above and a collection of distinct complex numbers $z_1, \dots, z_n$,
we may
construct one more vector space of the same dimension as
$\Sing\,L_{\bs \La}[\La^{(\infty)}]$ and an algebra of commuting linear
operators acting on that new space.

\medskip
Namely, write $\La^{(i)}=(\La^{(i)}_1,\dots,\La^{(i)}_N)$,
$i=1,\dots,n,\infty$, with
$\La^{(i)}_1 \geq \dots \geq \La^{(i)}_{N-1}\geq \La^{(i)}_N$ being
non-negative integers. Consider the vector space $\C_d[x]$ of
polynomials in $x$ of degree not greater than $d$, where $d$ is
a natural number big enough
with respect to $n$ and $N$. Define $n+1$ Schubert cycles
$C_{z_1,\La^{(1)}}, \dots, C_{z_n,\La^{(n)}},
C_{\infty,\La^{(\infty)}}$ in the Grassmannian of all $N$-dimensional
subspaces of $\C_d[x]$ as follows. For $i=1,\dots,n$,
the cycle $C_{z_i,\La^{(i)}}$ is the closure
of the set of all $N$-dimensional subspaces $V\subset \C_d[x]$ having a basis
$f_1,\dots,f_N$ such that
$f_j(x) = (x-z_i)^{\La^{(i)}_{j}+N-j} + O((x-z_i)^{\La^{(i)}_{j}+N-j+1})$
for all $j$.
The cycle $C_{\infty,\La^{(\infty)}}$ is the closure
of the set of all $N$-dimensional subspaces $V\subset \C_d[x]$ having a basis
$f_1,\dots,f_N$ of polynomials of
degrees $\La^{(\infty)}_N, \La^{(\infty)}_{N-1}+1, \dots,
\La^{(i)}_{\infty}+N-1$, respectively.
Consider the intersection
of these cycles and the algebra $A_G$ of functions on this
intersection.

By Schubert calculus, the dimension of $A_G$, regarded as a vector space, equals
the dimension of the vector space $\Sing\,L_{\bs \La}[\La^{(\infty)}]$.
Multiplication in the algebra $A_G$ defines on the vector space $A_G$ the
commutative algebra of linear multiplication operators. The vector space $A_G$
with the commutative algebra of multiplication operators is our new object.

We conjecture that there exists a natural isomorphism of the vector spaces
$A_G\to\Sing\,L_{\bs \La}[\La^{(\infty)}]$ which induces an isomorphism of
the corresponding algebras --- the algebra of multiplication operators on $A_G$
and the Bethe algebra $A_L$ acting on $\Sing\,L_{\bs \La}[\La^{(\infty)}]$.

\medskip
Note that the Bethe algebra $A_L$ has linear algebraic nature (it is generated
by a finite set of relatively explicitly defined matrices) while the algebra
$A_G$ has geometric nature (it is the algebra of functions on the intersection
of several algebraic cycles). An isomorphism of $A_L$ and $A_G$ may allow us
to study one of the algebras in terms of the other.

For example, the intersection of Schubert cycles
$C_{z_1,\La^{(1)}}, \dots, C_{z_n,\La^{(n)}}, C_{\infty,\La^{(\infty)}}$ is not
transversal if and only if the algebra $A_G$ has nilpotent elements. Probably
it is easier to check the presence of such elements in $A_L$ than in $A_G$.

As another example, assume that all elements of the Bethe algebra
$A_L$ are diagonalizable.
In that case the algebra $A_G$ does not have nilpotent elements, hence
the intersection of the Schubert cycles is transversal. Returning back to
the Bethe algebra $A_L$ we may conclude that the spectrum
of  $A_L$ is simple.

\medskip
The main result of this paper is the construction of an isomorphism
of $A_L$ and $A_G$ for $N=2$.

\subsection{}
The paper has the following structure.

\medskip
In Section \ref{Two algebras} we define two algebras $A_M$ and $A_D$.
The algebra $A_M$ is the algebra generated by the Gaudin Hamiltonians acting
of the subspace $\snash$ of singular vectors of weight $\La^{(\infty)}$ in the tensor
product $M_{\bs \La} = M_{\La^{(1)}} \otimes \dots \otimes M_{\La^{(n)}}$
of Verma $\glt$-modules. Here $\La^{(i)} = (m_s,0)$ for $i=1,\dots,n$ and
$\La^{(\infty)} = (\sum_{s=1}^nm_s - l, l)$.

To define the algebra $A_D$ we consider the scheme $C_D$ of monic
linear second-order differential operators on $\C$ having regular singular
points at $z_1,\dots,z_n, \infty$, with exponents $0, m_i+1$ at $z_i$
for $i=1,\dots, n$, and exponents $-l,l-1-\sum_{s=1}^nm_s$ at infinity,
and also having a polynomial of degree $l$ in its kernel. Then we define
$A_D$ as the algebra of functions on $C_D$.

In Section \ref{Epimorphism psi-DM} we construct an algebra
epimorphism $\psi_{DM} : A_D \to A_M$.

\medskip

In Section \ref{On separation for slt}
we describe Sklyanin's separation of variables
for the $\glt$ Gaudin model and introduce the universal weight function.
The important result of Section
\ref {On separation for slt}
is Theorem \ref{thm on Bethe ansatz} on the Bethe ansatz method,
which describes the interaction of the three objects: algebras
$A_M$, $A_D$, and the universal weight function.

\medskip

In Section \ref{Multiplication Mult} we
consider the space $A_D^*$, dual to the vector space $A_D$,
and the algebra of linear operators on $A_D^*$ dual to the
multiplication operators on $A_D$. Using the universal weight
function we construct a linear map $\tau : A_D^* \to
\Sing\,\nash$. Theorem \ref{1st main thm}
says that $\tau$ is an isomorphism identifying the algebra
of operators on $A^*_D$ dual to multiplication operators and
the Bethe algebra $A_M$ acting on $\snash$.
Theorem \ref{1st main thm} is our first main result.

\medskip

In Section \ref{Grothendieck bilinear form on A-D} using the
Grothendieck bilinear form on $A_D$ we construct an isomorphism $\phi
: A_D \to A^*_D$. The isomorphism $\phi$ identifies the algebra of multiplication
operators on $A_D$ with the algebra of operators on $A^*_D$ dual to
multiplication operators.

\medskip

In Section \ref{Three more algebras} we introduce three more algebras $A_G$, $A_P$,
$A_L$.

The algebra $A_G$ is the algebra of functions on the intersection of
Schubert cycles $C_{z_1,\La^{(1)}}, \dots, C_{z_1,\La^{(n)}},
C_{\infty,\La^{(\infty)}}$ in the Grassmannian of two-dimensional subspaces
of
$\C_d[x]$.

To define the algebra $A_P$ we consider the scheme $C_P$ of monic
linear second-order differential operators on $\C$ having regular singular
points at $z_1,\dots,z_n, \infty$, with exponents $0, m_i+1$ at $z_i$
for $i=1,\dots, n$ and exponents $-l,l-1-\sum_{s=1}^nm_s$ at infinity,
and also having the kernel consisting of polynomials only.
Then the algebra
$A_P$ is the algebra of functions on $C_P$.

The algebra $A_M$ is the algebra generated by the Gaudin Hamiltonians acting
of the subspace $\Snash$ of singular vectors of weight $\La^{(\infty)}$ in
the tensor product
$L_{\bs \La} = L_{\La^{(1)}} \otimes \dots \otimes L_{\La^{(n)}}$
of polynomial irreducible finite-dimensional $\glN$-modules
with highest weights $\La^{(1)}, \dots, \La^{(n)}$, respectively.

\medskip

In Section \ref{New homomorphisms} we discuss interrelations of the
five algebras $A_D, A_M, A_G, A_P, A_L$. In particular, we have
a natural isomorphism $\psi_{GP} :A_G \to A_P$.

In Section \ref{New homomorphisms}
we construct
a linear map $\zeta : A_P \to \Snash$. Using our first main result we
show in Theorem \ref{second main thm} that $\zeta$ is an
isomorphism identifying the algebra of multiplication operators on
$A_P$ and the Bethe algebra $A_L$ acting on $\Snash$.
Theorem \ref{second main thm} is our
second main result.

In Section \ref{Equations with polynomial solutions only} using the
Shapovalov form on $\Snash$ and the isomorphism $\zeta$ we construct a
linear map $\theta : A^*_P \to \Snash$. In Theorem \ref {useful cor}
we show that $\theta$ is an isomorphism identifying the algebra on
$A_P^*$ of operators dual to multiplication operators and the Bethe
algebra $A_L$ acting on $\Snash$. This is our third main result.

As an application of the third main result we prove
the following statement, see Corollary \ref{Cor 2 of 2nd main thm}.

\medskip
\noindent

{\it
If a two-dimensional vector space $V$
belongs to the intersection of the Schubert
cycles $C_{z_1,\La^{(1)}}, \dots, C_{z_1,\La^{(n)}},
C_{\infty,\La^{(\infty)}}$ and if
$d^2/dx^2 + a(x) d/dx + b(x)$ is the differential operator
annihilating $V$, then there exists a nonzero eigenvector $v\in \Snash$
of the Bethe algebra $A_L$ with eigenvalues given by the functions
$a(x)$ and $b(x)$.}

\medskip
Note that the converse statement follows from Corollaries 12.2.1 and 12.2.2 in
\cite{MTV3}, see Sections~\ref{Cor 2 of 2nd main thm NEW} and
\ref{Cor 2 of 2nd main thm}.

\medskip
In Appendix we discuss the relations between the Grothendieck residue on
$A_D$, the Shapovalov form on $\Snash$ and the homomorphism
$A_D \to \snash \to \Snash$.

\subsection{}
We thank P.\,Belkale and F.\,Sottile for useful discussions.

\section{Two algebras}
\label{Two algebras}

\subsection{Algebra $A_M$}
\label{Algebra B}

\subsubsection{}
\label{subsec gl2}
Let $\glt$ be the
complex Lie algebra of $2\times 2$-matrices with standard generators
$e_{ab}, a,b=1,2$.
Let $\h\subset \glt$ be the Cartan subalgebra of diagonal matrices,
$\h^*$ the dual space,
$(\,,\,)$ the standard scalar product on $\h^*$,\
$\epsilon_1,\epsilon_{2}\in \h^*$ the standard orthonormal basis,
$\al = \epsilon_1-\epsilon_{2}$ the simple root.

Let $\bs \La = (\La^{(1)},\dots, \La^{(n)})$ be a collection of
$\glt$-weights, where $\La^{(s)} = m_s\epsilon_1$ with $m_s\in \C$.

Let $l $ be a nonnegative integer. Define the $\glt$-weight\
$\La^{(\infty)} = \sum_{s=1}^n \,\La^{(s)} - l\,\alpha$.

\medskip
The pair $\bs \La$, $l$ is called {\it separating} if \
$\sum_{s=1}^n m_s - 2l +1 + i\, \neq \, 0$
for all $i = 1, \dots, l$.

\subsubsection{}
Let $\bs z = (z_1,\dots,z_n)$ be a collection of distinct complex numbers.
Let
\be
M_{\bs \La}
\ = \ M_{\La^{(1)}}\otimes \dots \otimes M_{\La^{(n)}}
\ee
be the tensor product of Verma $\glt$-modules with
highest weights $\La^{(1)}, \dots , \La^{(n)}$, respectively.
Denote by $\Sing M_{\bs\La}[\La^{(\infty)}]$ the subspace of
$M_{\bs\La}$ of singular vectors of weight $\La^{(\infty)}$,
\be
\Sing\,M_{\bs \La}[\La^{(\infty)}]\ =\
\{\,
v \in M_{\bs\La}\ |\ e_{12}v = 0,\ e_{22}v = l v \, \}\ .
\ee
Consider the differential operator
\be
\D_{M_{\bs \La}}\ =\
\left(\frac{d}{dx} -
\sum_{s=1}^n \frac {e_{11}^{(s)}}{x-z_s}\right)
\left(\frac{d}{dx} -
\sum_{s=1}^n \frac {e_{22}^{(s)}}{x-z_s}\right)
-
\left(\sum_{s=1}^n \frac {e_{21}^{(s)}}{x-z_s}\right)
\left(\sum_{s=1}^n \frac {e_{12}^{(s)}}{x-z_s}\right) .
\ee
The differential operator acts on $M_{\bs \La}$-valued functions in $x$
and is called {\it the universal differential operator} associated with
$M_{\bs \La}$ and $\bs z$, \cite{T}, \cite{MTV1}, \cite{MTV3}.
We have
\beq
\D_{M_{\bs \La}}\ =\ \frac{d^{2}}{dx^{2}} \ -\
\sum_{s=1}^n \,\frac{m_s}{x-z_s}\,\frac{d}{dx}\ +\
\sum_{s=1}^n\,\frac{\Ht_s}{x-z_s}\
\eeq
where $\Ht_1, \dots, \Ht_n \in \End\,(M_{\bs \La})$,
\beq
\label{form for Gaudin}
\Ht_s\ = \ \sum_{r\ne s}\ \frac 1
{z_s-z_r}\
(\,m_s m_r - \Omega_{s,r} \,)\
\qquad
{\rm and}
\qquad
\Omega_{s,r}\ =\ \sum_{i,j=1}^2 e_{ij}^{(s)}\otimes e_{ji}^{(r)}\ .
\eeq
We have $\Ht_1+\dots + \Ht_n = 0$.

The operators
$\Ht_1, \dots, \Ht_n$ are called {\it the Gaudin Hamiltonians}
associated with $M_{\bs \La}$ and $\bs z$.
The Gaudin Hamiltonians have the following properties:
\begin{enumerate}
\item[(i)]
The Gaudin Hamiltonians commute:\
$[\Ht_i,\Ht_j]=0$
for all $i,j$.

\item[(ii)]
The Gaudin Hamiltonians commute
with the $\glt$-action
on $M_{\bs \La}$:
$[\Ht_i,x]=0$ for all $i$ and $x\in U(\glt)$.
\end{enumerate}
In particular, the Gaudin Hamiltonians preserve
the subspace $\snash \subset M_{\bs \La}$.

Restricting $\D_{M_{\bs \La}}$ to the subspace of $\snash$-valued
functions we obtain the differential operator
\beq
\D_{\Sing M_{\bs \La}}
\ =\ \frac{d^{2}}{dx^{2}} \ -\
\sum_{s=1}^n \,\frac{m_s}{x-z_s}\,\frac{d}{dx}\ +\
\sum_{s=1}^n\,\frac{H_s}{x-z_s}\
\eeq
where $H_s=\Ht_s|_{\snash}\,\in \End\,(\snash)$.

The operator
$\D_{\Sing M_{\bs \La}}$ will be called {\it the universal
differential operator} associated with
$\snash$ and $\bs z$. The operators
$H_1,\dots,H_n$ will be called
{\it the Gaudin Hamiltonians}
associated with $\snash$ and $\bs z$.

\goodbreak
The commutative unital subalgebra
of $\End\,(\snash)$ generated by the Gaudin Hamiltonians $H_1,\dots,H_n$
will be called {\it the Bethe algebra} associated with $\snash$ and $\bs z$
and denoted by $A_M$.

\subsubsection{}
Introduce the operators $G_0,\dots, G_{n-2}$ by the formula
\be
\sum_{s=1}^n\,\frac{H_s}{x-z_s}\ =\
\frac { G_{0}x^{n-2} + \dots + G_{n-2}}{ (x-z_1)\dots (x-z_n)}\ .
\ee
Then
\ $G_0\, =\, l \,(\,\sum_{s=1}^n \,m_s\, +\,1\, - \, l)\,$.

\subsubsection{}
\label{dim if fixed}
{\bf Lemma.}\
{\it
Assume that the pair $\bs \La, l$\ is separating. Then}
\bea
&&
\dim\,\Sing\,M_{\bs \La}
\bigl[\,{\textstyle \sum_{s=1}^n \,\La^{(s)} \,-\, l\,\alpha}\>\bigr]\ =
\\
&&
\phantom{aaaaa}
\dim\,M_{\bs \La}
\bigl[\,{\textstyle \sum_{s=1}^n \,\La^{(s)} \,-\, l\,\alpha}\>\bigr]\ -\
\dim\,M_{\bs \La}
\bigl[\,{\textstyle \sum_{s=1}^n \,\La^{(s)} \,-\,(l-1)\,\alpha}\>\bigr]\ .
\eea

\begin{proof}
The map \ $
e_{12}e_{21} :
M_{\bs \La}\bigl[\,\sum_{s=1}^n \,\La^{(s)} - (l-1)\,\alpha\>\bigr]
\to
M_{\bs \La}\bigl[\,\sum_{s=1}^n \,\La^{(s)} - (l-1)\,\alpha\>\bigr]
$\
is an isomorphism of vector spaces since the pair
$\bs \La, l$ is separating.
The fact that $e_{12}e_{21}$ is an isomorphism implies the lemma.
\end{proof}

\subsubsection{}
\label{thm exist of polyn}
{\bf Theorem.}
{\it
Assume that the pair $\bs \La, l$\ is separating. Then for any\\
$v_0\in \snash$ there exist unique $v_1,\dots,v_l \in\snash$
such that the function
\be
v(x) \ =\ v_0\,x^l \, +\,v_1\,x^{l-1}\, +\, \dots \,+\,v_l
\ee
is a solution of the differential equation $\D_{\Sing M_{\bs \La}} v(x)\,=\,0$.
}

\begin{proof}
If all weights $\La^{(1)},\dots,\La^{(n)}$ are dominant integral, then
the theorem holds by Theorem~12.1.3 from \cite{MTV3}.
By Lemma \ref{dim if fixed} the dimension of $\snash$ does not depend on
$\bs\La$ if the pair $\bs\La,l$ is separating. Hence the theorem holds
for all separating $\bs\La,l$.
\end{proof}

\subsection{Algebra $A_D$}
\subsubsection{}
\label{subsub h's}
Denote $\bs a = (a_1,\dots,a_l)$ and $\bs h = (h_1,\dots,h_n)$.
Consider the space $\C^{l+n}$ with coordinates $\bs a, \bs h$.
Denote by $D$ the set of all points $\T\in \C^{l+n}$ whose coordinates
satisfy the equations $q_{-1}(\bs h) = 0, \ q_0(\bs h) = 0$, where
\be
\label{defining eqns for D}
q_{-1}(\bs h)\ = \ \sum^n _{s=1}\, h_s \, ,
\qquad
q_0(\bs h)\ = \ \sum^n _{s=1}\, z_s h_s\, - \,
l\,( \sum^n _{s=1} m_s + 1 - l)\ .
\ee
The set $D$ is an affine space of dimension $l+n-2$.

\subsubsection{}
Denote by $\Dh$ the following polynomial
differential operator in $x$ depending on parameters $\bs h$,
\beq
\label{D}
\Dh\ =\
\left(\prod_{s=1}^n\,(x-z_s)\right)\!\!
\left(\frac{d^{2}}{dx^{2}} -
\sum_{s=1}^n \frac{m_s}{x-z_s}\frac{d}{dx}
+ \sum_{s=1}^n \frac{h_s}{x-z_s} \right)\ .
\eeq
If $\bs p \in D$, then the singular points of
$\Dhp$ are $z_1,\dots,z_n,\infty$ and
the singular points are regular.
For $s = 1, \dots , n$, the exponents of $\Dhp$ at $z_s$ are
$0, m_s+1$.
The exponents of $\Dhp$ at $\infty$ are
$-l, l-1- \sum_{s=1}^nm_s$.

\subsubsection{}
\label{subsec p(x,a)}
Denote by $\qa$ the following polynomial in $x$ depending on
parameters $\bs a$,
\be
\qa \ =
\ x^l + a_1x^{l-1}+ \dots + a_l \ .
\ee
If $\bs h$ satisfies equations $q_{-1}(\bs h) = 0$
and $q_{0}(\bs h) = 0$, then
the polynomial $\Dh(\qa)$ is a polynomial in $x$ of degree $l+n-3$,
\be
\Dh(\qa)\ =\ q_1(\bs a,\bs h)\,x^{l+n-3}\ +\ \dots \ +\
q_{l+n-2}(\bs a,\bs h)
\ .
\ee
The coefficients $q_i (\bs a,\bs h)$ are functions linear
in $\bs a$ and linear in $\bs h$.

Denote by $I_D$ the ideal in
$\C[\bs a, \bs h]$ generated by
polynomials $q_{-1}, q_0, q_1,\dots,q_{l+n-2}$.
The ideal $I_D$ defines a scheme $C_D\subset D$. Then
\be
A_D\ =\ \C[\bs a, \bs h]/ I_D\
\ee
is the algebra of functions on $C_D$.

The scheme $C_D$ is the scheme of points $\T\in D$ such that the
differential equation
$\Dhp u(x)=0$ has a polynomial solution $p(x,\bs a(\T))$.

\subsubsection{}

The scheme $C_D$ and the algebra $A_D$
depend on the choice of distinct numbers $\bs z = (z_1,\dots,z_n)$:
\ ${}$ $C_D = C_D(\bs z)$,\ $A_D = A_D(\bs z)$.

\subsubsection{}
\label{thm const wrt z}
{\bf Theorem.}
{\it
Assume that the pair $\bs \La$, $l$ is separating.
Then
the dimension of $A_D(\bs z)$, considered as a
vector space, is finite and does not depend on the choice of distinct
numbers $z_1,\dots,z_n$.
}

\begin{proof}
It suffices to prove two facts:
\begin{enumerate}
\item[(i)]
For any $\bs z$ with distinct coordinates
there are no algebraic curves lying in $C_D(\bs z)$.

\item[(ii)]
Let a sequence $\bs z^{(i)}$, $i=1,2,\dots{}$, tend to a finite limit
$\bs z = (z_1,\dots,z_n)$ with distinct $z_1,\dots, z_n$. Let
$\T^{(i)} \in C_D(\bs z^{(i)})$, $i = 1, 2, \dots\ $,
be a sequence of points. Then all coordinates
$(\bs a(\T^{(i)}), \bs h(\T^{(i)})$ remain bounded as $i$ tends to infinity.
\end{enumerate}
We prove (i), the proof of (ii) is similar.

For a point $\T$ in $C_D(\bs z)$, the operator
$\D_{\bs h(\T)}$ has the form
\be
B_0(x)\frac{d^2}{dx^2} + B_1(x)\frac d{dx} + B_2(x,\T)
\ee
where the polynomials $B_0,B_1,B_2$ in $x$
are of degree $n,n-1,n-2$, respectively, the top degree coefficients
of the polynomials $B_0,B_1,B_2$
are equal to $1, -\sum_{s=1}^nm_s$, $l (\sum_{s=1}^n m_s + 1 - l)$, respectively,
and the polynomials $B_0, B_1$ do not depend on $\T$.

Assume that (i) is not true. Then there exists a sequence of points
$\T^{(i)} \in C_D(\bs z)$, $i=1,2,\dots{}$, which tends to infinity
as $i$ tends to infinity.

Then it is easy to see that $\bs h(\T^{(i)})$ cannot tend to infinity since
it would contradict to the fact that $\D_{\bs h(\T^{(i)})} ( p(x,\bs
a(\T^{(i)}))) = 0$.

Now choosing a subsequence we may assume that
$\bs h(\T^{(i)})$ has finite limit as $i$ tends to infinity.

If $\bs h(\T^{(i)})$ has finite limit as $i$ tends to infinity,
then $\bs a(\T^{(i)})$ cannot tend to infinity since it would
mean that the limiting differential equation has a polynomial solution of degree
less than $l$ and this is impossible.

This reasoning implies
that $\T^{(i)} \in C_D(\bs z)$ cannot tend to infinity. Thus we get contradiction
and statement (i) is proved.
\end{proof}

\subsection{Second description of $A_D$}
\subsubsection{}
\label{thm-conjecture 2}
{\bf Theorem.}
{\it
Assume that the pair $\bs \La, l$ is separating.
Assume that $\bs h$ satisfies equations $q_{-1}(\bs h) = 0$
and $q_{0}(\bs h) = 0$.
Consider the system
\beq
\label{system q}
q_i(\bs a, \bs h)\ =\ 0\ ,
\qquad
i = 1, \dots , l\ ,
\eeq
as a system of linear equations with respect to $a_1,\dots,a_l$.
Then this system has a unique solution
\ $a_i = a_i(\bs h)$,\ $i = 1, \dots, l$,\
where $a_i(\bs h)$ are polynomials in $\bs h$. }
\hfill
$\square$

\begin{proof}
Theorem \ref{thm-conjecture 2} follows from the fact that
\be
q_i(\bs a,\bs h) \ =\
i\, ( \sum_{s=1}^n m_s - 2l + i + 1 )\, a_i \ +\ \sum_{j=1}^{i-1} \,
q_{ij}(\bs h)\, a_j
\ee
for $i = 1, \dots , l$. Here $q_{ij}$ are some linear functions of $\bs h$.
The coefficient of $a_i$
does not vanish because the pair $\bs\La, l$ is separating.
\end{proof}

\subsubsection{}
\label{sec 2nd description}
Denote by $I'_D$ the ideal in
$\C[\bs h]$ generated by $n$ polynomials $q_{-1}, q_0$,
$q_j(\bs a(\bs h),\bs h)$, $j= l+1,\dots, l+n-2$. Then
\be
A_D \ \cong \ \C[\bs h]/ I'_D\ .
\ee

\subsection{Third description of $A_D$}

\subsubsection{}
Assume that $h_1,\dots,h_n$ satisfy equations $q_{-1}(\bs h)= 0,\
q_{0}(\bs h)= 0$.
Then
\be
\sum_{s=1}^n \frac {h_s}{x-z_s}\ =\
\frac {g(x)}{(x-z_1)\dots (x-z_n)}\ ,
\vvgood
\ee
where
\be
g(x)\ =\ l\,(\sum_{s=1}^n m_s + 1 - l)\, x^{n-2} +
g_1(\bs h)x^{n-3}+ g_2(\bs h)x^{n-2}+ \dots +
g_{n-2}(\bs h)\
\ee
for suitable $g_1(\bs h),\dots,g_{n-2}(\bs h)$ which are
linear functions in $\bs h$.

\subsubsection{}
\label{lemma-conjecture 3}
{\bf Lemma.}
{\it
Let $c_1,\dots,c_{n-2}$ be arbitrary numbers.
Consider the system of $n$ linear equations
\be
\sum^n _{s=1}\, h_s\, =\, 0 \, ,
\qquad
\sum^n _{s=1}\, z_s h_s\, = \,
l\,( \sum^n _{s=1} m_s + 1 - l)\ ,
\ee
\be
g_i(\bs h) = c_i
\qquad
i=1,\dots,n-2\ ,
\ee
with respect to $h_1,\dots,h_n$.
Then this system has a unique solution.}
\hfill
$\square$

\medskip
This lemma is the standard fact from the theory of simple fractions.

\subsubsection{}
Let $\bs g = (g_0,\dots,g_{n-2})$ be a tuple of numbers and
\vvn.4>
\be
g(x)\ =\
g_0x^{n-2} + g_1x^{n-3} + \dots + g_{n-2}\ .
\ee
The expression
\be
(\prod_{s=1}^n(x-z_s))
( \frac {d^2}{dx^2}p(x,\bs a)
- \sum_{i=1}^n \frac {m_i}{x-z_i}\frac{d}{dx}p(x,\bs a)) +
g(x)p(x,\bs a) \ =\ 0\ .
\ee
is a polynomial in $x$ of degree $l+n-2$,
\vvn.4>
\be
\hat q_0(\bs a,\bs g)\,x^{l+n-2}\ +\
\hat q_1(\bs a,\bs g)\,x^{l+n-3}\ +
\ \dots \ + \hat q_{l+n-2}(\bs a,\bs g)\ ,
\vv.4>
\ee
where $\hat q_0(\bs a,\bs g)\,=\,g_0-l\,(\sum_{s=1}^n m_s + 1 - l)$.

\subsubsection{}
\label{lemma-conjecture 4}
{\bf Lemma.}
{\it The system of equations
\be
\hat q_i(\bs a,\bs g)\ =\ 0\ ,
\qquad
i=0,\dots, n-2\ ,
\ee
determines $g_0,\dots,g_{n-2}$ uniquely as polynomials in $\bs a$.}
\hfill
$\square$

\begin{proof}
The equation $\hat q_0(\bs a,\bs g)\,=\,0$ gives
$g_0\,=\,l\,(\sum_{s=1}^n m_s + 1 - l)$.
Now Lemma \ref{lemma-conjecture 4} follows from the fact that
\be
\hat q_i(\bs a,\bs g)\ =\ g_i\ +\ \sum_{j=1}^{i-1}\, \hat q_{ij}(\bs a) g_j
\ee
for $i=1,\dots,n-2$. Here $\hat q_{ij}$ are some linear functions of $\bs a$.
\end{proof}

\goodbreak
\subsubsection{}
Combining Lemmas \ref{lemma-conjecture 3}
and
\ref{lemma-conjecture 4}, we obtain polynomial functions
$h_i = h_i(\bs a)$, $i=1,\dots, n$.

Denote by $I_D''$ the ideal in $\C[\bs a]$ generated by $l$
polynomials $q_j(\bs a,\bs h(\bs a))$, $j= n-1,\dots,
l+n-2$. Then
\be
A_D \ \cong \ \C[\bs a]/ I_D''\ .
\ee

\subsection{ Epimorphism $\psi_{DM} : A_D \to A_M$}
\label{Epimorphism psi-DM}

Let $h_1,\dots,h_n$ be the functions on $D$, introduced in
Section \ref{subsub h's}, and $H_{1}, \dots,H_n$
the Gaudin Hamiltonians.

\subsubsection{}
\label{thm D to B}
{\bf Theorem.}
{\it
Assume that the pair $\bs \La, l$\ is separating. Then the
assignment \\
$h_{s}\ \mapsto\ H_{s}$, \ $s = 1, \dots , n$,
determines an algebra
epimorphism\
$\psi_{DM} : A_{D} \to A_M$.
}

\begin{proof}
The equations defining the scheme $C_D$ are the equations of
existence of a polynomial solution
$\qa$ of degree $l$ to the polynomial
differential equation
$\D_{\bs h} u(x)\,=\,0$. By Theorem \ref{thm exist of polyn},
the defining equations for $C_D$ are satisfied by the coefficients of
the universal differential operator $\D_{\Sing M_{\bs \La}}$.
\end{proof}

\section{Separation of variables}
\label{On separation for slt}

\subsection{Holomorphic representation}

The tensor product
$M_{\bs \La} = M_{\La^{(1)}} \otimes \dots \otimes M_{\La^{(n)}}$
of Verma $\glt$-modules is identified with the space of polynomials
$\C[x^{(1)},\dots,x^{(n)}]$ by the linear map
\be
e_{21}^{j^1}v_{\La^{(1)}} \otimes \dots \otimes
e_{21}^{j^n}v_{\La^{(n)}}\ \mapsto\ (x^{(1)})^{j^1}\dots (x^{(n)})^{j^n}\,,
\ee
where $v_{\La^{(s)}}$ is the generating vector of $M_{\La^{(s)}}$.
Then the $\glt$-action on $\C[x^{(1)},\dots,x^{(n)}]$
is given by the differential operators,
\be
e_{12}^{(s)} = - x^{(s)}\p_{x^{(s)}}^2 + m_s \partial_{x^{(s)}} \ ,
\qquad
e_{21}^{(s)} = x^{(s)} \ ,
\ee
\be
e_{11}^{(s)} = - 2x^{(s)}\partial_{x^{(s)}} + m_s \ ,
\qquad
e_{22}^{(s)} = 0 \ ,
\ee
where $\partial_{x^{(s)}}$ denotes the derivative with respect to
$ x^{(s)}$.

\subsection{Change of variables}
\label{change of var slt}

Make the change of variables from
$x^{(1)},\ldots,x^{(n)}$ to
$ u$, $y^{(1)}$,
\dots, $y^{(n-1)}$
using the relation
\be
\label{change}
\sum^n_{s=1}\ \frac{x^{(s)}}{t-z_s}\ =\
u\,
\frac{ \prod^{n-1}_{k=1}\, (t - y^{(k)})}
{\prod^n_{s=1}\, (t-z_s)}\ ,
\ee
where $t$ is an indeterminate. This relation defines $u, y^{(1)},
\ldots,y^{(n-1)}$
uniquely up to permutation of $y^{(1)}, \dots, y^{(n-1)}$
unless
$u = \sum_{s=1}^n x^{(s)} =0$.
The map
$
( u, y^{(1)},\dots,y^{(n-1)})\mapsto (x^{(1)},
\dots,x^{(n)})$
is an unramified covering on the complement to the union of diagonals
$y^{(i)}=y^{(j)}$, $i\neq j$, and the hyperplane $u=0$.

\subsection{Sklyanin's theorem}

Consider the operators $\Ht_1,\dots, \Ht_n$
defined by formula \Ref{form for Gaudin}.
Introduce the operators
\be
K_i (\Ht) \ =\ \sum_{s=1}^n\ \frac 1 {y^{(i)} - z_s}\ \Ht_s\ ,
\qquad i=1,\dots,n-1\ .
\ee

\subsubsection{}
\label{Sklyanin thm}
{\bf Theorem \cite{Sk}.}
{\it In variables $u, y^{(1)},\dots,y^{(n-1)}$, we have}
\be
K_i(\Ht) \ =\ -\,\p^2_{y^{(i)}}\ +
\ \sum_{s=1}^n \, \frac{m_s}{y^{(i)}-z_s}\, \p_{y^{(i)}}\ ,
\qquad
i = 1,\dots, n-1\ .
\ee

\subsection{Universal weight function}
\label{Universal weight function}
The weight subspace
$M_{\bs \La}[\La^{(\infty)}] \subset M_{\bs \La}$
is identified with
the subspace of
$\C[x^{(1)}, \dots, $ $x^{(n)}]$ of homogeneous
polynomials of degree $l$.

We consider the associated $\nash$-valued
universal weight function
\be
\prod_{j=1}^{l}\,(\,\prod_{i=1}^n (t_{j}-z_i)
\sum_{s=1}^n\,
\frac {x^{(s)}} { t_{j}-z_s}\,)\
\ee
of variables
$x^{(1)},\dots,x^{(n)}$, $t_1,\dots,t_{l}$.
In variables $u, y^{(1)},\dots, y^{(n-1)}$,
$t_1,\dots,t_{l}$,
the universal weight function takes the form\
$
(-1)^{ln}\,
u^{l}\, \prod_{j=1}^{n-1}\,
p(y^{(j)}),
$
\ where
$p(x) = \prod_{i=1}^{l}\,(x-t_i)$. If we denote by $-a_1, a_2, \dots,
(-1)^l a_l$ the elementary symmetric functions of $t_1,\dots,t_l$, then
$p(x) = p(x,\bs a)$ in notation of Section
\ref{subsec p(x,a)},
and the universal weight function takes the form
\be
\omega(u, \bs y, \bs a)\ = \
(-1)^{ln}\,
u^{l}\, \prod_{j=1}^{n-1}\, p(y^{(j)}, \bs a)\ ,
\ee
with $\bs y = (y^{(1)}, \dots, y^{(n-1)})$.

The trivial but important property of the universal weight function is
given by the following lemma.

\subsubsection{}
\label{lemma omega is nonzero}
{\bf Lemma.}
{\it For every $\T \in D$, the vector $\omega(u, \bs y, \bs a(\T))$ is
a nonzero vector of $\nash$}.
\hfill
$\square$

\medskip
Denote by $\omega_D$ the projection of the universal weight function
$\omega(u, \bs y, \bs a)$ to $M_{\bs \La} \otimes A_D$.

\subsubsection{}
\label{thm on Bethe ansatz}
{\bf Theorem.}
{\it For $s = 1, \dots , n$, we have
\beq
\label{formula H equals h}
\Ht_s \,\omega_D\ =\ h_s \,\omega_D\
\eeq
in \ $M_{\bs \La}\otimes A_D$. Moreover, we have }
\beq
\label{formula omega is sing}
\omega_D\ \in\ \Sing M_{\bs \La}[\La^{(\infty)}]\otimes A_D\ .
\eeq

\begin{proof}
First we prove formula \Ref{formula H equals h}.
Let $\C(u,\bs y)$ be the algebra of rational functions in $u,\bs y$.
For $i=1,\dots,n-1$, introduce
\be
K_i (\bs h) \ =\ \sum_{s=1}^n\ \frac {h_s} {y^{(i)} - z_s}\
\in \C(u,\bs y)\otimes A_D \ .
\ee
We claim that
\beq
\label{neW}
K_i (\Ht)\,\omega_D\ =\ K_i (\bs h)\,\omega_D
\eeq
in $\C(u,\bs y)\otimes A_D$. Indeed,
\vvn.2>
\begin{align*}
K_i & (\Ht)\, \omega(u, \bs y, \bs a) \
= \ (K_i (\bs h)
+ K_i (\Ht) - K_i (\bs h))\,\omega(u, \bs y, \bs a)\ =
\ K_i (\bs h) \,\omega(u, \bs y, \bs a)\ +{}
\\[4pt]
& (-1)^{ln}\,
u^{l}\left[\left(-\,\p^2_{y^{(i)}}\ +
\ \sum_{s=1}^n \, \frac{m_s}{y^{(i)}-z_s}\, \p_{y^{(i)}}\ - \
\sum_{s=1}^n\ \frac 1 {y^{(i)} - z_s}\ h_s\right) p(y^{(i)},\bs a)\right]
\prod_{j\neq i} p(y^{(j)},\bs a) .
\end{align*}
Clearly, the last term has zero projection to
$\C (u,\bs y)\otimes A_D$ and we get formula \Ref{neW}.

Having formula \Ref{neW}, let us show that $\Ht_s \omega_D \,= \,h_s \omega_D$
in $\C[u, \bs y]\otimes A_D$. For that introduce two
$\C[u, \bs y]\otimes A_D$-valued functions in a new variable $x$:
\vvn.2>
\be
F_1(x)\ =\ \sum_{s=1}^n \frac{\Ht_s\omega_D}{x-z_s}\ ,
\qquad F_2(x)\ =\ \sum_{s=1}^n \frac{h_s\omega_D}{x-z_s}\ ,
\ee
and show that the functions are equal.

Each of the functions is the ratio of a polynomial in $x$ of degree $n-2$
and the polynomial $(x-z_1)\dots (x-z_n)$. To check that the two functions
are equal it is enough to check
that $F_1(x)=F_2(x)$ for $x=y^{(i)}$, $i=1,\dots,n-1$, but this follows from
formula \Ref{neW}. Hence formula \Ref{formula H equals h} is proved.

Formula \Ref{formula omega is sing} follows from
formula \Ref{formula H equals h}.
Indeed, by formula \Ref{form for Gaudin} we have
$\sum_{s=1}^n z_s \Ht_s =
\sum_{s=1}^n \sum_{r=1}^{s-1} \,(m_s m_r - \Omega_{s,r})$.
This implies that $\sum_{s=1}^n z_s \Ht_s$
acts on the weight subspace $M_{\bs\La}[\La^{(\infty)}]$ as
the operator
$l ( \sum_{s=1}^n m_s + 1 - l ) - E_{21} E_{12}$,
where $E_{ij} = \sum_{s=1}^n e^{(s)}_{ij}$. Since\
$\sum_{s=1}^n\, z_s h_s = l ( \sum_{s=1}^n m_s + 1 - l )$,
formula \Ref{formula H equals h} allows us to conclude that\
$E_{21} E_{12}\, \omega_D = 0$.
The operator $E_{21}$ is injective, in variables $u, y^{(1)},\dots,y^{(n-1)}$
it is the operator of multiplication by $u$.
Therefore, $E_{12} \,\omega_D\, =\, 0$.
\end{proof}

\section{Multiplication in $A_D$ and Bethe algebra $A_M$}
\label{Multiplication Mult}

\subsection{Multiplication in $A_D$}

By Theorem \ref{thm const wrt z}, the scheme $C_D$ considered as a set is
finite, and the algebra $A_D$ is the direct sum of local algebras
corresponding to points $\T$ of the set $C_D$,
\be
A_D\ =\ \oplus_{\T}\ A_{\T, D}\ .
\vv.3>
\ee
The local algebra $\AD$ may be defined as the quotient of the algebra of germs
at $\T$ of holomorphic functions in $\bs a, \bs h$ modulo the ideal $I_{\T,D}$
generated by all functions $q_{-1},\dots,q_{l+n-2}$.
The local algebra $\AD$ contains the maximal ideal $\mg_\T$
generated by germs which are zero at $\T$.

For $f\in A_D$, denote by $L_f$ the linear operator
$A_D \to A_D, \ g\mapsto fg$,
\ of multiplication by $f$. Consider the dual space
\be
A_D^*\ =\ \oplus_\T\, \AD^*
\vv.2>
\ee
and the dual operators $L_f^* : A_D^* \to A_D^*$.
Every summand $\AD^*$ contains the distinguished one-dimensional subspace
$\mg^\T$ which is the annihilator of $\mg_\T$.

\subsubsection{}
\label{lemma on dual operators}
{\bf Lemma.}
{\it
\begin{enumerate}
\item[(i)]
For any point $\T$ of the scheme
$C_D$ considered as a set and any $f\in A_D$, we have
$L_f^* (\mg^\T) \subset \mg^\T$.

\item[(ii)] For any point $\T$ of the scheme
$C_D$ considered as a set, if $W\subset \AD^*$ is a nonzero vector subspace
invariant with respect to all
operators $L_f^*$, $f\in A_D$, then $W$ contains $\mg^\T$.
\end{enumerate}
}

\begin{proof} For any $f \in \mg_\T$ we have
$L_f^*(\mg^\T) = 0$. This gives part (i).

To prove part (ii) we
consider the filtration of $\AD$ by powers of the
maximal ideal,
\be
\AD \supset \mg_\T \supset \mg_\T^2\supset\dots\supset \{0\}\ .
\ee
We consider a linear basis $\{f_{a,b}\}$ of $\AD$, $a=0,1,\dots{}$,
$b=1,2,\dots{}$, which agrees with this filtration. Namely, we assume that
for every $i$, the subset of all vectors $f_{a,b}$ with $a\geq i$ is a basis
of $\mg^i_\T$\,.

Since dim $\AT/\mg_\T = 1$, there is only one basis vector with $a=0$ and
we also assume that this vector $f_{0,1}$ is the image of $1$ in $\AD$\,.

Let $\{ f^{a,b}\}$ denote the dual basis of $\AD^*$. Then
the vector $f^{0,1}$ generates $\mg^\T$.

Let $w = \sum_{a,b} c_{a,b} f^{a,b}$ be a nonzero vector in $W$. Let
$a_0$ be the maximum value of $a$ such that there exists $b$ with a nonzero
$c_{a,b}$. Let $b_0$ be such that $c_{a_0,b_0}$ is nonzero.
Then it is easy to see that $L^*_{f_{a_0,b_0}} w\,=\,c_{a_0,b_0} f^{0,1}$.
Hence $W$ contains $\mg^\T$.
\end{proof}

\subsection{Linear map $\tau : A_D^* \to \Sing\,\nash$}

Let $f_1,\dots, f_\mu$ be a basis of $A_{D}$
considered as a vector space over $\C$. Write
\beq
\label{basis in M}
\omega_D\ =\ \sum_i v_i
\otimes f_i
\qquad
{\rm with} \qquad v_i \in\snash \ .
\eeq
Denote by $V \subset\snash$ the vector subspace spanned by
$v_1,\dots,v_\mu$.
Define the linear map
\beq
\label{map psi}
\tau\ :\ A_{D}^*\ \to\ \snash\ ,
\qquad
g\ \mapsto\ g(\omega_D) = \sum_i\ g(f_i)\,v_i\ .
\eeq
Clearly, $V$ is the image of $\tau$.

\subsubsection{}
\label{lemma april 1}
{\bf Lemma.}
{\it Let $\T$ be a point of $C_D$ considered as a set. Let
$\omega(u, \bs y, \bs a(\T)) \in \nash$ be the value of the universal weight
function at $\T$. Then the vector $\omega(u, \bs y, \bs a(\T))$ belongs to
the image of $\tau$.}
\hfill
$\square$

\subsubsection{}
\label{lemma 1}
{\bf Lemma.}
{\it Assume that the pair $\bs \La, l$\ is separating. Then for
any $f\in A_D$ and $g\in A_D^*$,
we have}\
$\tau (L^*_f (g)) = \psi_{DM}(f) (\tau(g))$.

\medskip
In other words, the map $\tau$ intertwines the action of the algebra
of multiplication operators $L^*_f$ on $A_{D}^*$ and the action on
the Bethe algebra on $\snash$.

\goodbreak
\begin{proof}
The algebra $A_D$ is generated by $h_1,\dots,h_n$. It is enough to prove
that for any $s$ we have $\tau (L^*_{h_s} (g))\, =\, H_s (\tau(g))$.
But $\tau(L^*_{h_s} (g)) = \sum_i g(h_s f_i) v_i =
g( \sum_i \,v_i \otimes h_sf_i ) = g(\sum_i \,H_s v_i \otimes f_i ) =
H_s (\tau(g))$.
\end{proof}

\subsubsection{}
\label{corollary 1}
{\bf Corollary.}
{\it
The vector subspace $V\subset\snash$
is invariant with respect to the action of the Bethe algebra $A_M$
and the kernel of $\tau$ is a subspace of $A_D^*$, invariant with respect to
multiplication operators $L^*_f,\,f\in A_D$.
}

\subsection{First main theorem}
\subsubsection{}
\label{1st main thm}
{\bf Theorem.}
{\it
Assume that the pair $\bs \La, l$\ is separating. Then the
image of $\tau$ is $\Sing\,\nash$ and the kernel of $\tau$ is zero.}

\subsubsection{}
\label{1st main cor}
{\bf Corollary.}
{\it
The map $\tau$ identifies the action of operators $L_f^*$, $f\in A_D$,
on $A_D^*$ and the action of the Bethe algebra on $\Sing\,\nash$.
Hence the epimorphism $\psi_{DM}\, : \, A_D\, \to\, A_M$\ is an isomorphism.}

\begin{proof}[Proof of Theorem \ref{1st main thm}]
Let $d\,=\,\dim\,\snash$. Theorem 9.16 in \cite{RV} says that for generic $\bs z$
there exists $d$ points $\T_1,\dots,\T_d$ in $C_D$ such that the vectors
$\omega(u, \bs y, \bs a(\T_1))$, \dots , $\omega(u, \bs y, \bs a(\T_d))$ form a basis in
$\snash$. Hence, $\tau$ is an epimorphism for generic $\bs z$ by
Lemma \ref{lemma april 1}. By Theorem \ref{thm const wrt z} and
Lemma \ref{dim if fixed}
dimensions of $A_D$ and $\snash$ do not depend on $\bs z$. Hence
$\dim\,A_D \geq \dim\,\snash$.
Therefore, to prove Theorem \ref{1st main thm} it is
enough to prove that $\tau$ has zero kernel.

Denote the kernel of $\tau$ by $K$.
Let $A_D = \oplus_\T \AD$ be the decomposition into the direct sum of local
algebras.
Since $K$ is invariant with respect to multiplication operators, we have
\ $K \,= \,\oplus_\T\, K\cap \AD^*$\ and for every $\T$ the vector subspace
$K\cap \AD^*$ is invariant with respect to multiplication operators.
By Lemma \ref{lemma on dual operators}, if $K\cap \AD^*$ is nonzero, then
$K\cap \AD^*$ contains the one-dimensional subspace $\mg^\T$.

Let $\{f_{a,b}\}$ be the basis of $\AD$ constructed in the proof of Lemma
\ref{lemma on dual operators} and let $\{f^{a,b}\}$ be the dual basis of
$\AD^*$. Then the vector $f^{0,1}$ generates $\mg^\T$. By definition of $\tau$,
the vector $\tau(f^{0,1})$ is equal to the value of the universal weight
function at $\T$. By Lemma \ref{lemma omega is nonzero}, this value is nonzero
and that contradicts to the assumption that $f^{0,1} \in K$.
\end{proof}

\subsection{Grothendieck bilinear form on $A_D$}
\label{Grothendieck bilinear form on A-D}

Realize the algebra $A_D$ as $\C[\bs h]/I'_D$, \
where $I'_D$ is the ideal generated by $n$ polynomials
$q_{-1}, q_0$, $q_j(\bs a(\bs h),\bs h)$, $j= l+1,\dots, l+n-2$,
see Section \ref{sec 2nd description}.

Let $\rho: A_{D} \to \C$, be the Grothendieck residue,
\be
f \ \mapsto \ \frac 1{(2\pi i)^n}\,\Res_{C_D}
\frac{ f}{q_{-1}(\bs h) q_0(\bs h)\prod_{j=l+1}^{l+n-2}\,
q_j(\bs a(\bs h),\bs h)}\ .
\vv.2>
\ee
Let $(\,,\,)_D$ be the Grothendieck symmetric
bilinear form on $A_{D}$ defined by the rule
\beq
( f,\, g)_{D}\ =\ \rho (f g)\ .
\eeq
The Grothendieck bilinear form is non-degenerate.

The form $(\,,\,)_D$ determines a linear isomorphism $\phi:A_D \to A_D^*$,
$f \mapsto (f,\,\cdot)_D$.

\subsubsection{}
\label{lem on intertw}
{\bf Lemma.}
{\it
The
isomorphism $\phi$ intertwines the operators $L_f$ and $L^*_f$
for any $f\in A_D$. }

\begin{proof} For $g\in A_D$ we have $\phi (L_f(g))=
\phi (fg) = (fg,\cdot)_D = (g, f\cdot)_D = L_f^*((g,\cdot)_D) =
L_f^*\phi(g)$.
\end{proof}

\subsubsection{}
\label{intersting cor}
{\bf Corollary.} {\it
Assume that the pair $\bs \La, l$\ is separating. Then the
composition $\tau\phi\,:\, A_D \to \snash$
is a linear isomorphism which intertwines the algebra of multiplication
operators on $A_D$ and the action of the Bethe algebra $A_M$ on
$\snash$.
}

\section{Three more algebras}
\label{Three more algebras}

\subsection{New conditions on $\bs \La,\,l$}
\label{sec new conditiong}
In the remainder of the paper we assume that
$\bs \La = (\La^{(1)},\dots, \La^{(n)})$ is a collection of dominant integral
$\glt$-weights,
\beq
\label{weights formula}
\La^{(s)}\ = \ m_s\,\epsilon_1\ ,
\qquad m_s \in \Z_{\geq 0}\ ,
\qquad
s = 1,\dots, n\ .
\eeq
We assume that $l \in \Z_{\geq 0}$ is such that the weight
$\La^{(\infty)}\, =\, \sum_{s=1}^n \La^{(s)}\, -\, l \alpha$
is dominant integral. Hence the pair $\bs \La, l$ is separating.

\subsection{Algebra $A_P$}
Denote $\tilde l = \sum_{s=1}^n m_s\,+1 - l$. We have $\tilde l > l$.
Denote
\be
\tilde{\bs a} = (\tilde a_1, \dots, \tilde a_{\tilde l - l -1},
\tilde a_{\tilde l - l +1}, \dots, \tilde a_{\tilde l})\ .
\ee
Consider space $\C^{\tilde l + l + n - 1}$
with coordinates $\tilde{\bs a}, \bs a, \bs h$,
cf.~Section \ref{subsub h's}.

Denote by $\qat$ the following polynomial in $x$ depending on parameters
$\tilde{\bs a}$,
\be
\qat \ =
\ x^{\tilde l} + \tilde {a}_1x^{\tilde l-1}+ \dots +
\tilde a_{\tilde l - l -1} x^{l+1} +
\tilde a_{\tilde l - l +1}x^{l-1}+ \dots +\tilde a_{\tilde l}\ .
\vvgood
\ee
If $\bs h$ satisfies the equations $q_{-1}(\bs h) = 0$
and $q_{0}(\bs h) = 0$, then the polynomial $\Dh(\qat)$ is a polynomial in $x$
of degree $\tilde l+n-3$,
\be
\Dh(\qat)\ =\ \tilde {q}_1(\tilde{\bs a},\bs h)\,x^{\tilde l+n-3}
\ +\ \dots \ +\ \tilde q_{\tilde l+n-2}(\tilde{\bs a},\bs h)\ .
\vv.2>
\ee
The coefficients $\tilde{q}_i (\tilde{\bs a},\bs h)$ are functions linear
in $\tilde{\bs a}$ and linear in $\bs h$.

Recall that if $\qa = x^l + a_1x^{l-1}+ \dots + a_l$ and $\bs h$ satisfies
equations $q_{-1}(\bs h) = 0$ and $q_{0}(\bs h) = 0$, then the polynomial
$\Dh(\qa)$ is a polynomial in $x$ of degree $l+n-3$,
\vvn.2>
\be
\Dh(\qa)\ =\ q_1(\bs a,\bs h)\,x^{l+n-3}\ +\ \dots \ +\
q_{l+n-2}(\bs a,\bs h)\ .
\vv.2>
\ee
Denote by $I_P$ the ideal in $\C[\tilde{\bs a}, \bs a, \bs h]$ generated
by polynomials $q_{-1}, q_0, q_1,\dots,q_{l+n-2}$,\\
$\tilde{q}_1,\dots,\tilde q_{\tilde l+n-2}$.

The ideal $I_P$ defines a scheme $C_P\subset \C^{\tilde l + l + n - 1}$.
The algebra
\be
A_P\ =\ \C[\tilde{\bs a}, \bs a, \bs h]/ I_P\
\ee
is the algebra of functions on $C_P$.

The scheme $C_P$ is the scheme of points $\T\in \C^{\tilde l + l + n - 1}$
such that the differential equation $\Dhp u(x)=0$ has two polynomial solutions
$\tilde {p}(x,\tilde{\bs a}(\T))$ and $p(x,\bs a(\T))$.

\subsection{Algebra $A_{G}$}
\label{Algebra AG}

Let $d$ be a sufficiently large natural number and
$\C_d[x]$ the vector subspace in $\C[x]$ of polynomials
of degree not greater than $d$.
Let $G$ be the Grassmannian of all two-dimensional vector subspaces in
$\C_d[x]$.
Let $\bs z = (z_1,\dots,z_n)$ be distinct complex numbers.

For $s = 1, \dots , n$, denote by $C_{z_s, \La^{(s)}} \subset G$ the Schubert
cycle associated with the point $z_s\in \C$ and weight $\La^{(s)}$.
The cycle $C_{z_s, \La^{(s)}}$ is the closure of the set
$C^o_{z_s, \La^{(s)}} \subset G$ of all two-dimensional subspaces
$V \subset \C_d[x]$ having a basis $f_1,f_2$ such that
\be
f_1(z_s) = 1 \qquad {\rm and} \qquad
f_2(x) = (x-z_s)^{m_s+1} + O((x-z_s)^{m_s+2})\ .
\ee

Denote by $C_{\infty, \La^{(\infty)}} \subset G$ the Schubert cycle
associated with the point $\infty$ and weight $\La^{(\infty)}$.
$C_{\infty, \La^{(\infty)}}$ is the closure of the set
$C^o_{\infty, \La^{(\infty)}} \subset G$ of all two-dimensional subspaces
$V \subset \C_d[x]$ having a basis $f_1,f_2$ such that
$\deg\, f_1 = l$ and $\deg\,f_2 = \tilde l$.

Consider the intersection
\be
C_G\ =\ C_{\infty, \La^{(\infty)}} \ \cap\ (\,\cap_{i=1}^n\,
C_{z_i, \La^{(i)}}\,)\ .
\ee
Denote by $A_G$ the algebra of functions on $C_G$.

It is known from Schubert calculus that $\dim \,A_G$ is finite and
does not depend on $\bs z$ with distinct coordinates.

\subsubsection{}
It is easy to see that
\be
C_G\ =\ C^o_{\infty, \La^{(\infty)}} \ \cap\ (\,\cap_{i=1}^n\,
C^o_{z_i, \La^{(i)}}\,)\ .
\ee

\subsubsection{}
We shall use the following presentation of the algebra $A_G$.

Consider space $\C^{\tilde l + l - 1}$ with coordinates $\tilde{\bs a}, \bs a$.
A point $\T \in\C^{\tilde l + l - 1}$ will be called admissible if for every
$s=1,\dots,n$ at least one of the numbers $\tilde {p}(z_s,\tilde{\bs a}(\T))$,
$p(z_s,\bs a(\T))$ is not zero. The set of all admissible points form a Zariski
open subset $U \subset \C^{\tilde l + l - 1}$.

\medskip

For polynomials $f,g \in \C[x]$ denote by $\Wr (f,g)$ the Wronskian $f'g-fg'$,
where ${}'$ denotes $d/dx$. The Wronskian of $\qat$ and $\qa$ has the form
\vv.2>
\be
\Wr\,(\qat, \qa) \ =\
(\tilde l - l) x^{\tilde l + l -1} +
w_1(\tilde{\bs a}, \bs a)x^{\tilde l + l -2}
+ \dots + w_{\tilde l + l -1}(\tilde{\bs a}, \bs a)\
\vv.2>
\ee
for suitable polynomials $w_1,\dots, w_{\tilde l + l -1}$
in variables $\tilde{\bs a}, \bs a$.

Let us write
\be
(\tilde l - l)\prod_{s=1}^n (x-z_s)^{m_s}\ =\
(\tilde l - l) x^{\tilde l + l -1} + c_1x^{\tilde l + l -2}
+ \dots + c_{\tilde l + l -1}\
\ee
for suitable numbers $c_1,\dots,c_{\tilde l + l -1}$.

Let $A_U$ be the algebra of regular functions on the set $U$
of all admissible points. Denote by $I_G \subset A_U$
the ideal generated by $\tilde l + l - 1$
polynomials $w_1-c_1,\dots, w_{\tilde l + l -1}-c_{\tilde l + l-1}$.
Then
\be
A_G \ =\ A_U/ I_G\ .
\ee
In this presentation of $A_G$ the scheme $C_G$ is the scheme of points
$\T \in U$ such that the Wronskian of $\tilde {p}(x,\tilde{\bs a}(\T))$ and
$p(x,\bs a(\T))$ is equal to $(\tilde l - l)\prod_{s=1}^n (x-z_s)^{m_s}$.

\subsection{Algebra $A_L$}

Let
\be
L_{\bs \La}\ = \ L_{\La^{(1)}}\otimes \dots \otimes L_{\La^{(n)}}
\vv.2>
\ee
be the tensor product of irreducible $\glt$-modules with
highest weights $\La^{(1)}, \dots , \La^{(n)}$, respectively.
Denote by $\Sing L_{\bs\La}[\La^{(\infty)}]$ the subspace of
$L_{\bs\La}$ of singular vectors of weight $\La^{(\infty)}$.

Let $S$ denote the tensor Shapovalov form on $\snash$, induced from
the tensor product of the Shapovalov forms on the factors of
$M_{\bs \La} = M_{\La^{(1)}} \otimes \dots \otimes M_{\La^{(n)}}$.

The Shapovalov form
determines the linear epimorphism
\vv.2>
\be
\sh \ :\ \snash \ \to\ \Snash\ .
\vv.2>
\ee
The Bethe algebra $A_M$ preserves the kernel of $\sh$ and induces a commutative
subalgebra $A_L$ in $\End\,(\Snash)$ called the Bethe algebra on $\Snash$.

Denote by $\psi_{ML} : A_M \to A_L$ the corresponding epimorphism.

\subsubsection{}
Denote by
\be
\D_{L}\ =\ \frac{d^{2}}{dx^{2}}
\ -\
\sum_{s=1}^n \,\frac{m_s}{x-z_s}\,\frac{d}{dx}\ +\
\sum_{s=1}^n\,\frac{\psi_{ML}(H_s)}{x-z_s}\
\ee
the universal differential operator associated with the subspace
$\Snash$ and collection $\bs z$.

\subsubsection{}
\label{thm exist of polyn in Snash}
{\bf Theorem.}
{\it
Assume that the pair $\bs \La, l$\ satisfies conditions of
Section \ref{sec new conditiong}.
Then for any
$v_0\in \Snash$ there exist $v_1,\dots,v_{\tilde l} \in \Snash$ such that
the function
\be
v(x) \ =\ v_0\,x^{\tilde l} \, +\,v_1\,x^{\tilde l-1}\, +\, \dots \,+\,
v_{\tilde l}
\ee
is a solution of the differential equation $\D_{L} v(x)\,=\,0$.
}

\medskip
This theorem is a particular case of Theorem 12.3 in \cite{MTV3}.

\section{Four more homomorphisms}
\label{New homomorphisms}

\subsection{Isomorphism $\psi_{GP}: A_G \to A_P$}
A point $\T$ of $C_P$ defines the differential equation
$\D_{\bs h(\T)} u(x) = 0$ and two solutions
$\tilde {p}(x,\tilde{\bs a}(\T))$ and $p(x,\bs a(\T))$. We have
\be
\Wr\,(\tilde {p}(x,\tilde{\bs a}(\T)),p(x,\bs a(\T)))\ =
\ (\tilde l - l)\prod_{s=1}^n (x-z_s)^{m_s}\ .
\ee
Hence, the pair $\tilde {p}(x,\tilde{\bs a}(\T))$, $p(x,\bs a(\T))$
defines a point of $C_G$.

This construction defines a homomorphism of algebras $\psi_{GP}: A_G \to A_P$.

\subsubsection{}
\label{thm 2 [BMV]}
{\bf Theorem.}
{\it
The homomorphism $\psi_{GP}$ is an isomorphism. }

\begin{proof} We construct the inverse homomorphism as follows. Let $\bs v$ be a point
of $C_G$. Consider the following differential equation with respect to a function
$u(x)$,
\be
\det
\left(
\begin{array}{ccc}
u'' & u' & u
\\
\tilde {p}(x,\tilde{\bs a}(\bs v))'' &
\tilde {p}(x,\tilde{\bs a}(\bs v))'
& \tilde {p}(x,\tilde{\bs a}(\bs v))
\\
p(x,\bs a(\bs v))'' & p(x,\bs a(\bs v))' & p(x,\bs a(\bs v))
\end{array}
\right) = \ 0\ .
\ee
Let us write this differential equation as
$B_0(x) u'' + B_1(x) u' + B_2(x) u = 0$. Here
\be
B_0(x)\ =\ \Wr(\tilde {p}(x,\tilde{\bs a}(\bs v)), p(x,\bs a(\bs v)))\ =\
(\tilde l - l)\prod_{s=1}^n (x-z_s)^{m_s}\ .
\ee
It is easy to see that each of the polynomials $B_1, B_2$ is divisible by
the polynomial
\be
B(x)\ =\ (\tilde l - l) \prod_{s=1}^n (x-z_s)^{m_s-1}\ .
\ee
Introduce the differential operator
\be
\D_{\bs v}\ =\ b_0(x) \frac{d^2}{dx^2}
+ b_1(x) \frac{d}{ dx} + b_2(x)\ =
\ \frac {1}{B(x)} \left( B_0(x)\frac{d^2}{ dx^2}
+ B_1(x)\frac{d}{ dx} + B_2(x)\right)\ .
\ee
Then
\be
b_0(x)\ = \ \prod_{s=1}^n \,(x-z_s)\ ,\qquad
b_1(x)\ = \ \prod_{s=1}^n \,(x-z_s)\left(
\sum_{s=1}^n\frac{-\>m_s}{x-z_s}\right)\ ,
\ee
and $b_2(x)$ is a polynomial of degree $n-2$, whose leading coefficient is
$\tilde ll$.

The triple, consisting of the differential operator $\D_{\bs v}$ and two
polynomials $\tilde {p}(x,\tilde{\bs a}(\bs v))$ and $p(x,\bs a(\bs v))$,
determines a point of $C_P$, thus defining the inverse homomorphism
$A_P \to A_G$.
\end{proof}

\subsubsection{}
\label{cor of thm 2 [BMV]}
{\bf Corollary.}
{\it
The dimension of the algebra $A_P$ is finite and does not depend on $\bs z$
with distinct coordinates.}

\medskip
Indeed, $\dim\,A_P = \dim\,A_G$ and $\dim \,A_G$ is finite and
does not depend on $\bs z$ with distinct coordinates.

\subsubsection{}
It is known from Schubert calculus that $\dim\,A_G = \dim\, \Snash$.

\subsection{Epimorphism $\psi_{DP} : A_D \to A_P$}

A point $\T$ of $C_P$ determines the differential equation $\Dhp\,u(x)\,=\,0$
and two solutions $\tilde {p}(x,\tilde{\bs a}(\T))$ and $p(x,\bs a(\T))$.
Then the pair, consisting of the differential equation $\Dhp\,u(x)\,=\,0$
and one of the solutions $p(x,\bs a(\T))$ determines a point of $C_D$. This
correspondence defines a natural algebra epimorphism $\psi_{DP} : A_D \to A_P$.

\subsection{Linear map $\xi : A_D \to \Snash$}

Denote by $\xi : A_D \to \Snash$ the composition of linear maps
\be
A_D\
\stackrel{\phi}{\longrightarrow}
\ A_D^*\
\stackrel{\tau}{\longrightarrow}
\ \snash\
\stackrel{\sh}{\longrightarrow}
\ \Snash\ .
\ee
By Theorem \ref{1st main thm}, \ $\xi$ is a linear epimorphism.

Denote by $\psi_{DL} : A_D \to A_L$ the algebra epimorphism defined as the
composition $\psi_{ML}\psi_{DM}$.

\subsubsection{}
\label{lEmma on intertw }
{\bf Lemma.}
{\it
The linear map $\xi$ intertwines the
action of the multiplication operators $L_f,\,f\in A_D$,\,
on $A_D$ and
the action of the Bethe algebra $A_L$ on $\Snash$, i.e. for any
$f,g\in A_D$ we have\ $\xi(L_f(g)) \,=\, \psi_{DL}(f)(\xi(g))$.
}

\medskip
The lemma follows from
Corollary \ref{intersting cor}.

\subsubsection{}
\label{tarasov's lem}
{\bf Lemma.}
{\it
The kernel of $\xi$ coincides with the kernel of $\psi_{DL}$.
}

\begin{proof}
If $\psi_{DL}(f)=0$, then
$
\xi(f) = \xi(L_f(1)) = \psi_{DL}(f)(\xi(1)) = 0.
$
On the other hand, if $\xi(f)=0$, then for any $g\in A_D$ we have
$
\psi_{DL}(f)(\xi(g)) = \xi(L_f(g)) = \xi(fg) = \xi(L_g(f))
= \psi_{DL}(g)(\xi(f)) = 0.
$
Since $\xi$ is an epimorphism, this means that $\psi_{DL}(f)=0$.
\end{proof}

\subsubsection{}
\label{2nd main lem}
{\bf Lemma.}
{\it The kernel of\, $\xi$ coincides with the kernel
of $\psi_{DP}$.
}

\begin{proof}
By Schubert calculus
$\dim\,\Snash\,=\,\dim\, A_G$. Hence it suffices to show that the
kernel of $\xi$ contains the kernel of $\psi_{DP}$. But this follows
from Theorems \ref{thm exist of polyn} and \ref{thm exist of polyn in
Snash}.

Indeed the defining relations in $A_P = A_D/(\ker\,\psi_{DP})$ are the
conditions on the operator $\Dh$ to have two linearly independent polynomials
in the kernel. Theorems
\ref{thm exist of polyn} and \ref{thm exist of polyn in Snash}
guarantee these relations for elements of the Bethe algebra
$A_L$.
Hence, the kernel of $\psi_{DL}$ contains the kernel of $\psi_{DP}$.
By Lemma \ref{tarasov's lem}, the kernel of $\xi$ coincides with
the kernel of $\psi_{DL}$. Therefore,
the kernel of $\xi$ contains the kernel of $\psi_{DP}$.
\end{proof}

\subsubsection{}
\label{Cor of 2nd main lem}
{\bf Corollary.}
{\it
Since the algebra epimorphisms $\psi_{DP}$ and $\psi_{DL}$ have
the same kernels, the algebras $A_P$ and $A_L$ are isomorphic, and hence
by Theorem~\ref{thm 2 [BMV]} the algebras $A_G$ and $A_L$ are isomorphic.
\hfill
$\square$
}

\medskip
\subsection{Second main theorem}

Denote by $\psi_{PL} : A_P \to A_L$ the isomorphism induced by $\psi_{DL}$ and
$\psi_{DP}$. The previous lemmas imply the following theorem.

\subsubsection{}
\label{second main thm}
{\bf Theorem.}
{\it
The linear map $\xi$ induces a linear isomorphism
\be
\zeta\ :\ A_P\ \to\ \Snash
\ee
which
intertwines the multiplication operators $L_f,\,f\in A_P$,\, on $A_P$ and
the action of the Bethe algebra $A_L$ on $\Snash$, i.e. for any $f,g\in A_P$
we have $\zeta(L_f(g)) \,=\, \psi_{PL}(f)(\zeta(g))$.
}
\hfill
$\square$

\subsubsection{}
\label{Cor 3 of 2nd main thm}
{\bf Corollary. }
{\it
If every operator $ f \in A_L$ is diagonalizable, then
the algebra $A_L$ has simple spectrum and
all of the points of the intersection of Schubert cycles
\be
C_G\ =\ C_{\infty, \La^{(\infty)}} \ \cap\ (\,\cap_{i=1}^n\,
C_{z_i, \La^{(i)}}\,)\
\ee
are of multiplicity one.
}

\begin{proof}[Proof of Corollary]
The algebras $A_L$, $A_P$ and $A_G$ are all isomorphic.
We have $A_P = \oplus_{\T} \,A_{\T,P}$ where the sum is over the points of the
scheme $C_P$ considered as a set and $A_{\T,P}$ is the local algebra associated
with a point $\T$. The algebra $A_{\T,P}$ has nonzero nilpotent elements
if $\dim\,A_{\T,P}>1$. If every element $f\in A_P$ is diagonalizable,
then the algebra $A_P$ is the direct sum of one-dimensional local algebras.
Hence $A_P$ has simple spectrum as well as the algebras $A_L$ and $A_G$.
\end{proof}

\subsubsection{}
\label{cor on shapiro}
Corollary \ref{Cor 3 of 2nd main thm} has the following application.

\medskip
\noindent
{\bf Corollary}\, \cite{EGSV}.
{\it If $\,z_1, \dots, z_n$ are real and distinct, then
all of the points of the intersection of Schubert cycles
\be
C_G\ =\ C_{\infty, \La^{(\infty)}} \ \cap\ (\,\cap_{i=1}^n\,
C_{z_i, \La^{(i)}}\,)\
\ee
are of multiplicity one.}

\medskip
\noindent
\begin{proof}
If $z_1,\dots, z_n$ are real and distinct, then by Corollary 3.5 in \cite{MTV2}
all elements of the Bethe algebra $A_L$ are diagonalizable operators. Hence the
spectrum of $A_G$ is simple and all points of $C_G$ are of multiplicity one.
\end{proof}

This corollary is proved in \cite{EGSV} by a different method.

\section{Operators with polynomial kernel and Bethe algebra $A_L$}
\label{Equations with polynomial solutions only}

\subsection{ Linear isomorphism $\theta : A_P^* \to \Snash$}

Define the symmetric bilinear form
on $A_P$ by the formula
\be
(f,\,g)_P\ = \ S(\zeta(f), \,\zeta(g))\qquad
{\rm for\ all}\quad f,g \in A_P \ .
\ee
Recall that $S(\,,\,)$ denotes the Shapovalov form.

\subsubsection{}
{\bf Lemma.}
{\it
The form $(\,,\,)_P$ is non-degenerate.}

\medskip

The lemma follows from the fact that the Shapovalov form on $\Snash$ is
non-degenerate and the fact that $\zeta$ is an isomorphism.

\subsubsection{}
{\bf Lemma.}
{\it
We have $(fg,h)_P = (g,fh)_P$ for all $f,g,h \in A_P$.}
\hfill
$\square$

\medskip
The form $(\,,\,)_P$ defines a linear isomorphism $\pi:A_P \to A^*_P$,
$f\mapsto (f\,,\cdot)_P$.

\subsubsection{}
\label{COR}
{\bf Corollary.}
{\it
The map $\pi$ intertwines the multiplication operators $L_f,\,f\in A_P$,\,
on $A_P$ and the dual operators $L^*_f,\,f\in A_P$,\, on $A^*_P$.
}

\subsection{Third main theorem}
Summarizing Theorem \ref{second main thm}
and Corollary \ref{COR}
we obtain the following theorem.

\subsubsection{}
\label{useful cor}
{\bf Theorem.}
{\it
The composition $\theta=\zeta\pi^{-1}$ is a linear isomorphism from $A^*_P$ to
$\Snash$ which intertwines the multiplication operators $L^*_f,\,f\in A_P$,\,
on $A^*_P$ and the action of the Bethe algebra $A_L$ on $\Snash$,
i.e. for any $f\in A_P$ and $g\in A^*_P$ we have
$\theta(L^*_f(g))\,=\, \psi_{PL}(f)(\theta(g))$.
}
\hfill
$\square$

\subsubsection{}
\label{Cor 2 of 2nd main thm NEW}
Assume that $v\in \Snash$ is an eigenvector of the Bethe algebra $A_L$,
that is, $\psi_{ML}(H_s) v = \lambda_s v$ for suitable $\lambda_s\in\C$
and $s=1,\dots,n$. Then, by Corollaries~12.2.1 and~12.2.2 in~\cite{MTV3},
the differential operator
\be
\D\ =\ \frac{d^{2}}{dx^{2}} \ -\
\sum_{s=1}^n \,\frac{m_s}{x-z_s}\,\frac{d}{dx}\ +\
\sum_{s=1}^n\,\frac{\lambda_s}{x-z_s}\
\ee
has the following properties.
The operator $\D$ has regular singular points at $z_1,\dots,z_n,\infty$.
For $s = 1, \dots , n$, the exponents of $\D$ at $z_s$ are $0, m_s+1$.
The exponents of $\D$ at $\infty$ are $-l, l-1- \sum_{s=1}^nm_s$.
The kernel of $\D$ consists of polynomials only.
The following corollary of Theorem \ref{useful cor} gives the converse
statement.

\subsubsection{}
\label{Cor 2 of 2nd main thm}
{\bf Corollary of Theorem \ref{useful cor}.}
{\it
Let $\T \in \C^n$ be a point such that $q_{-1}(\bs h(\T)) = 0$,
$q_0(\bs h(\T)) = 0$, and all solutions of the differential equation
$\Dhp u(x) = 0$ are polynomials. Then there exists an eigenvector
$v \in \Snash$ of the action of the Bethe algebra $A_L$ such that
for every $s=1,\dots,n$ we have}
\be
\psi_{ML}(H_s) \,v \,=\, h_s(\T)\, v\ .
\ee

\begin{proof}[Proof of Corollary \ref{Cor 2 of 2nd main thm}]
Indeed, such $\T$ defines a linear function $\eta : A_P \to \C$,
$h_s \mapsto h_s(\T)$ for $s=1,\dots,n$. Moreover, $\eta(fg) = \eta(f)\eta(g)$
for all $f,g \in A_P$. Hence $\eta \in A_P^*$ is an eigenvector of
multiplication operators on $A_P^*$. By Theorem \ref{useful cor}
this eigenvector corresponds to an eigenvector $v\in\Snash$ of the action
of the Bethe algebra $A_L$ with eigenvalues prescribed
in Corollary~\ref{Cor 2 of 2nd main thm}.
\end{proof}

\subsubsection{}
\label{how to find}
Assume that $\T \in \C^n$ is a point satisfying the assumptions of
Corollary \ref{Cor 2 of 2nd main thm}. We describe how to find the eigenvector
$v\in\Snash$ indicated in Corollary \ref{Cor 2 of 2nd main thm}.

Let $f(x)$ be the
monic polynomial of degree $l$ which is a solution of the differential
equation $\Dhp w(x) = 0$. Consider the polynomial
\be
\omega(u, \bs y)
\ = \ u^l \prod_{j=1}^{n-1} f(y^{(j)})
\ee
as an element of $M_{\bs \La}$, see Section \ref{Universal weight function}.
By Theorem \ref{thm on Bethe ansatz} this vector lies in
$\snash$ and $\omega(u,\bs y) $ is an eigenvector of the Bethe algebra
$A_M$ with eigenvalues presecribed in Corollary \ref{Cor 2 of 2nd main thm}.
Consider the maximal subspace $V \subset \snash$ with three properties:
\ i)\, $V$ contains $\omega(u, \bs y)$,\
ii)\,
$V$ does not contain other eigenvectors of the Bethe algebra $A_M$,
\ iii)\,
$V$ is invariant with respect to the Bethe algebra $A_M$.
\ Let $\sh (V) \subset \Snash$ be the image of $V$ under the epimorphism
$\sh$. Then the subspace $ \sh(V)$ contains a unique one-dimensional subspace
of eigenvectors of the Bethe algebra $A_L$. Any
such an eigenvector may serve as
an eigenvector of the Bethe algebra $A_L$ indicated in
Corollary \ref{Cor 2 of 2nd main thm}.

\section{Appendix. Grothendieck and Shapovalov forms}

\subsection{Form $(\,,\,)_S$ on $A_D$}

Define the symmetric bilinear form on $A_D$ by
the formula
\be
(f,\,g)_S\ = \ S(\xi(f), \,\xi(g))
\qquad
{\rm for\ all}
\quad
f,g \in A_D \ ,
\ee
where $S(\,,\,)$ denotes the Shapovalov form.

\subsubsection{}
{\bf Lemma.}
{\it
The kernel of the bilinear form $(\,,\,)_S$ coincides with the kernel
of the linear map $\xi$. }

The lemma follows from the fact that the Shapovalov form on $\Snash$
is non-degenerate.

\subsubsection{}
{\bf Lemma.}
{\it
We have
$(fg, h)_S = (g, fh)_S$ for all $f,g,h \in A_D$.
}

\medskip
The lemma follows from Theorem \ref{1st main thm} and the fact that
the operators of the Bethe algebra are symmetric with respect to the
Shapovalov form, see, for example, \cite{RV} and \cite{MTV1}.

\subsubsection{}
\label{cor of lemma 2}
{\bf Corollary.}
{\it There exists $F \in A_D$ such that\
$(f,g)_S\, =\, (Ff,g)_D$ for all $f,g \in A_F$.
}

\subsubsection{}
\label{lemma 3}
{\bf Lemma.}
{\it
The kernel of the multiplication operator
$L_F : A_D \to A_D$ coincides with the kernel of $\xi$.
}

\medskip
The lemma follows from Theorem \ref{1st main thm} and the fact that
the kernel of $\sh$ is the kernel of the Shapovalov form
on $\snash$.

\medskip

The image of $L_F$ is the principal ideal $(F) \subset A_D$ generated by $F$.

\subsubsection{}
\label{cor of lem 3}
{\bf Corollary.}
{\it
The algebra of operators $L_f, f\in A_D$,
restricted to $(F)$ is isomorphic to the algebra $A_L$.}

\medskip

Denote $J = \{ f\in A_D\ |\ fg=0 \ {\rm for\ all}\ g\in\ker\psi_{DP} \}$.
The following lemma describes the ideal $(F)$ without using
the Shapovalov form.

\subsubsection{}
{\bf Lemma.}
{\it
We have $(F)=J$.
}

\begin{proof}
The inclusion $(F) \subset J$ follows from Lemmas \ref{lemma 3}
and
\ref{2nd main lem}.
On the other hand, since
$(\,,\,)_D$ is non-degenerate, we have
$\dim J$ $=$ $ \dim A_D$ $-$ $\dim\ker\psi_{DP}$.
By Lemma \ref{lemma 3}, $(F)$ has the same dimension and hence $(F) = J$.
\end{proof}

\end{document}